\documentclass[a4paper,twoside,reqno,11pt]{article}

\usepackage{amssymb,amsmath,amsthm}
\usepackage{graphicx}
\usepackage{times}
\usepackage{color}

\usepackage[round]{natbib}

\leftmargini=\baselineskip

\newtheoremstyle{localthm}
	{7pt} % space above
	{7pt} % space below
	{\sl} % Body font
	{} % Indent amount
	{\bf} % Theorem head font
	{{\rm.}} % Punctuation after theorem head
	{.7em} % Space after theorem head
	{} % Theorem head spec ?

\theoremstyle{localthm}
\newtheorem{thm}{Theorem}%[section]

\newtheorem{lem}[thm]{Lemma}

\newtheoremstyle{localrem}
	{5pt} % space above
	{5pt} % space below
	{\rm} % Body font
	{} % Indent amount
	{\bf} % Theorem head font
	{{\rm.}} % Punctuation after theorem head
	{.7em} % Space after theorem head
	{} % Theorem head spec ?

\theoremstyle{localrem}
\newtheorem{exmp}[thm]{Example}

\newcommand{\bs}{\boldsymbol}

\newcommand{\bbeta}{\bs{\beta}}
\newcommand{\beps}{\bs{\varepsilon}}
\newcommand{\bmu}{\bs{\mu}}
\newcommand{\ba}{\bs{a}}
\newcommand{\bb}{\bs{b}}
\newcommand{\bB}{\bs{B}}
\newcommand{\bH}{\bs{H}}
\newcommand{\T}{\bs{T}}
\newcommand{\bu}{\bs{u}}

\newcommand{\bv}{\bs{v}}

\newcommand{\x}{\bs{x}}
\newcommand{\X}{\bs{X}}
\newcommand{\y}{\bs{y}}
\newcommand{\Y}{\bs{Y}}
\newcommand{\z}{\bs{z}}
\newcommand{\Z}{\bs{Z}}
\newcommand{\R}{\mathbb{R}}

\newcommand{\Ex}{\operatorname{\mathbb{E}}}
\renewcommand{\Pr}{\operatorname{\mathbb{P}}}
\newcommand{\spann}{\operatorname{\mathrm{span}}}

\newcommand{\bbM}{\mathbb{M}}
\newcommand{\bbO}{\mathbb{O}}
\newcommand{\bbS}{\mathbb{S}}
\newcommand{\bbV}{\mathbb{V}}
\newcommand{\Frm}{\mathrm{F}}
\newcommand{\Haar}{\operatorname{\mathrm{Haar}}}

\begin{document}
\addtolength{\baselineskip}{0.35\baselineskip}
%==================================================================

\title{Connecting Model-Based and Model-Free Approaches\\
	to Linear Least Squares Regression}
\author{Lutz D\"umbgen\footnote{Work supported by Swiss National Science Foundation}
	\ and Laurie Davies\\
	(University of Bern and University of Duisburg-Essen)}
\date{December 2024}
\maketitle

\begin{abstract}
In a regression setting with a response vector and given regressor vectors, a typical question is to what extent the response is related to these regressors, specifically, how well it can be approximated by a linear combination of the latter. Classical methods for this question are based on statistical models for the conditional distribution of the response, given the regressors. In the present paper it is shown that various p-values resulting from this model-based approach have also a purely data-analytic, model-free interpretation. This finding is derived in a rather general context. In addition, we introduce equivalence regions, a reinterpretation of confidence regions in the model-free context.
\end{abstract}

%---------------------
\section{Introduction}
\label{sec:intro}
%---------------------

Statistical inference with general linear models is a well-established and indespensable tool for data analysis. The standard output of statistical software for linear models includes least squares estimators of parameters and their standard erros as well as p-values for various linear hypotheses. While the latter are based on certain model assumptions, linear models can also be viewed as tools for purely exploratory data analysis. In such a model-free context, one might wonder whether the p-values for, say, the relevance of certain covariates are still meaningful. The surprising answer is yes, these p-values do have a very precise and new interpretation.

To formulate a first result, suppose we observe a response vector $\y \in \R^n$ and $p < n$ linearly independent regressor vectors (regressors) $\x_1, \ldots, \x_p \in \R^n$. For a given integer $p_o \in \{0,\ldots,p-1\}$, we would like to know whether the least squares approximation of $\y$ by a linear combination $\hat{\y}$ of $\x_1,\ldots,\x_p$, that is,
\[
	\bigl\| \y - \hat{\y} \bigr\|^2 \
	= \ \min_{\bs{\beta} \in \R^p} \,
		\Bigl\| \y - \sum_{j=1}^p \beta_j \x_j \Bigr\|^2
\]
with the standard Euclidean norm $\|\cdot\|$, is substantially better than the restricted least squares approximation of $\y$ by a linear combination $\hat{\y}_o$ of $\x_1,\ldots,\x_{p_o}$ only, where $\hat{\y}_o := \bs{0}$ in case of $p_o = 0$. A classical, model-based answer is to assume that the $\x_j$ are fixed while
\begin{equation}
\label{eq:Gaussian.model}
	\y \ \sim \ N_n \Bigl( \sum_{j=1}^p \beta_j \x_j, \sigma^2 \bs{I} \Bigr)
\end{equation}
with an unknown parameter vector $\bs{\beta} \in \R^p$ and an unknown standard deviation $\sigma > 0$. Then an exact p-value of the null hypothesis that
\[
	\beta_j = 0 \ \ \text{for} \ p_o < j \le p
\]
is given by
\begin{equation}
\label{eq:p-value.model.F}
	1 - \mathrm{F}_{p-p_o,n-p} \Bigl(
		\frac{\|\hat{\y} - \hat{\y}_o\|^2/(p - p_o)}{\|\y - \hat{\y}\|^2/(n - p)} \Bigr) ,
\end{equation}
where $\mathrm{F}_{k,\ell}$ denotes the distribution function of Fisher's F distribution with $k$ and $\ell$ degrees of freedom. Since $\|\y - \hat{\y}_o\|^2 = \|\y - \hat{\y}\|^2 + \|\hat{\y} - \hat{\y}_o\|^2$, one can deduce from well-known connections between chi-squared, gamma, F and beta distributions that the p-value \eqref{eq:p-value.model.F} may be rewritten as
\begin{equation}
\label{eq:p-value.model.beta}
	\mathrm{B}_{(n-p)/2,(p - p_o)/2} \Bigl(
		\frac{\|\y - \hat{\y}\|^2}{\|\y - \hat{\y}_o\|^2} \Bigr) ,
\end{equation}
where $\mathrm{B}_{a,b}$ denotes the distribution function of the beta distribution $\mathrm{Beta}(a,b)$ with parameters $a, b > 0$.

Now let us view all observation vectors $\y$ and $\x_1,\ldots,\x_p$ as fixed. To measure to what extent $(\x_j)_{p_o < j \le p}$ contributes substantially to the least squares fit $\hat{\y}$, let $\hat{\y}^*$ be the least squares fit of $\y$ after replacing $(\x_j)_{p_o < j \le p}$ with a tuple $(\x_j^*)_{p_o < j \le p}$ of independent random vectors $\x_j^* \sim N_n(\bs{0},\bs{I})$. A precise measure of the relevance of $(\x_j)_{p_o < j \le p}$ is given by the probability that $\|\y - \hat{\y}^*\|^2$ is not larger than $\|\y - \hat{\y}\|^2$. The smaller this probability, the higher is the relevance. Interestingly, it can be computed exactly and coincides with the p-value \eqref{eq:p-value.model.beta}.

\begin{lem}
\label{lem:p-value.model-free.beta}
For arbitrary fixed, linearly independent vectors $\y, \x_1,\ldots,\x_p \in \R^n$ and stochastically independent standard Gaussian random vectors $\x_j^* \in \R^n$, $p_o < j \le p$,
\[
	\frac{\|\y - \hat{\y}^*\|^2}{\|\y - \hat{\y}_o\|^2} \
	\sim \ \mathrm{Beta} \bigl( (n-p)/2, (p-p_o)/2 \bigr) .
\]
In particular,
\[
	\Pr \bigl( \|\y - \hat{\y}^*\|^2 \le \|\y - \hat{\y}\|^2 \bigr) \
	= \ \mathrm{B}_{(n-p)/2,(p - p_o)/2} \Bigl(
		\frac{\|\y - \hat{\y}\|^2}{\|\y - \hat{\y}_o\|^2} \Bigr) .
\]
\end{lem}

This lemma is essentially a variant of a classical result about the angle between a random linear subspace and a fixed vector, see for instance Theorem~1.1 of \cite{Frankl_Maehara_1990}. A direct and self-contained proof will be given in Appendix~\ref{app:proof.F.lemma}. But Lemma~\ref{lem:p-value.model-free.beta} can be viewed as a special case of a more general connection between the model-based and model-free point of view which is elaborated in Section~\ref{sec:F-test.etc}. In particular, the classical model-based p-values do not require a Gaussian distribution of $\y$, given $(\x_j)_{1 \le j \le p}$, and for the model-free interpretation, the random tuple $(\x_j^*)_{p_o < j \le p}$ may have different distributions all of which lead to the p-value \eqref{eq:p-value.model.beta}. In Section~\ref{sec:Equivalence.regions} we discuss ``equivalence'' regions. In the model-based context, these are confidence regions for the unknown mean vector $\bmu = \Ex(\y)$. Under the model-free point of view, the interpretation of these regions is somewhat different. To illustrate the concept, we describe relatively simple equivalence regions for a sparse signal vector.

Some final comments and hints to additional work are given in Section~\ref{sec:Final.comments}. In particular, we explain how the considerations and results in the present paper are related to previous work about permutation tests in regression settings, a key reference being \cite{Freedman_Lane_1983} and a review of \cite{Winkler_etal_2014}.

Technical details and proofs are deferred to the appendix. Throughout this paper we use standard results from multivariate statistics and linear models as presented in standard textbooks, e.g.\ \cite{Mardia_Kent_Bibby_1979}, \cite{Eaton_1983} and \cite{Scheffe_1959}, without further reference.

%-----------------------------------------------
\section{The F test and other methods revisited}
\label{sec:F-test.etc}
%-----------------------------------------------

We consider arbitrary vectors $\y$ and $\x_1,\ldots,\x_p$ in $\R^n$. At first we discuss the question wether there is any association between $\y$ and $(\x_j)_{1 \le j \le p}$. In the introduction, this corresponds to $p_o = 0$. At the end of this section (Section~\ref{subsec:Composite.null.models}) we return to situations in which the contribution of $p_o \in \{1,\ldots,p-1\}$ regressors $\x_1,\ldots,\x_{p_o}$ is not questioned. This includes linear regression models with an intercept, accommodated by the trivial regressor $\bs{x}_1 = (1)_{i=1}^n$.

Concerning the regressors $\x_j$, suppose the raw data are given by a data matrix with $n$ rows
\[
	[y_i, \bs{w}_i^\top] \ = \ [y_i, w_{i,1}, \ldots, w_{i,d}] , \quad 1 \le i \le n ,
\]
containing the values of a response and $d$ covariates for each observation. If the covariates are numerical or $0$-$1$-valued, the usual multiple linear regression model would consider the regressors $\x_1 := (1)_{i=1}^n$ and $\x_j := (w_{i,j-1})_{i=1}^n$, $2 \le j \le d+1$. More complex models would also include the $\binom{d}{2}$ interaction vectors $\x_{j(a,b)} := (w_{i,a} w_{i,b})_{i=1}^n$, $1 \le a < b \le d$. In general, with arbitrary types of covariates, one could think of $\x_j = (f_j(\bs{w}_i))_{i=1}^n$ with given basis functions $f_1,\ldots,f_p$.

Let us introduce some notation. The unit sphere of $\R^n$ is denoted by $\bbS_n$, and $\bbO_n$ stands for the set of orthogonal matrices in $\R^{n\times n}$.

%------------------------------------
\subsection{The model-based approach}
%------------------------------------

We consider the regressors $\x_1,\ldots,\x_p$ as fixed and $\y$ as a random vector. In settings with random regressors, the subsequent considerations concern the conditional distribution of $\y$, given $\x_1,\ldots,\x_p$. For simplicity we assume throughout that the distribution of $\y$ is continuous, i.e.\ $\Pr(\y = \y_o) = 0$ for any fixed $\y_o \in \R^n$.

The null hypothesis of no relationship between $\y$ and the regressors $\x_1,\ldots,\x_p$ can be specified by describing a distribution of $\y$ which does not depend on the latter vectors (or any other fixed regressors). 

\smallskip

\noindent
$H_o$: \ The random vector $\y$ has a spherically symmetric distribution on $\R^n$. That means, its length $\|\y\|$ and direction $\|\y\|^{-1} \y$ are stochastically independent, where $\|\y\|^{-1} \y \sim \mathrm{Unif}(\bbS_n)$, the uniform distribution on the unit sphere $\bbS_n$.

\smallskip

It is well-known that this hypothesis $H_o$ encompasses the classical assumption that $\y \sim N_n(\bs{0}, \sigma^2 \bs{I})$ for some unknown $\sigma > 0$.

\paragraph{P-values for $H_o$.}
Let $S(\y) = S(\y,\x_1,\ldots,\x_p)$ be a test statistic such that high values indicate a potential violation of $H_o$. Then a p-value for $H_o$ is given by
\begin{equation}
\label{eq:pv.model-based.1}
	\pi(\y) \ := \ \Pr \bigl( S(\|\y\| \bu) \ge S(\y) \,\big|\, \y \bigr) ,
\end{equation}
where $\bu \sim \mathrm{Unif}(\bbS_n)$ is independent from $\y$. If $S(\y)$ is scale-invariant in the sense that
\begin{equation}
\label{eq:scale-invariance}
	S(c\bv) \ = \ S(\bv)
	\quad\text{for all} \ \bv \in \R^n \setminus \{\bs{0}\} \ \text{and} \ c > 0,
\end{equation}
one can write
\begin{equation}
\label{eq:pv.model-based.2}
	\pi(\y) \ = \ 1 - \Frm(S(\y)\,-) ,
\end{equation}
with the distribution function $\Frm$ of $S(\bu)$,
\[
	\Frm(x) \ := \ \Pr(S(\bu) \le x) .
\]
Here one could also consider a random vector $\z \sim N_n(\bs{0},\bs{I})$ instead of $\bu$.

\begin{exmp}[F test]
\label{ex:F}
If $\x_1,\ldots,\x_p$ are linearly independent with $p < n$, and if $S(\y)$ equals the F test statistic
\begin{equation}
	\label{eq:F.statistic}
	S(\y) \ := \ \frac{\|\hat{\y}\|^2/p}{\|\y - \hat{\y}\|^2/(n - p)} ,
\end{equation}
then scale-invariance of the latter implies that the p-value $\pi(\y)$ is given by the simplified formula \eqref{eq:pv.model-based.2}. Moreover, the distribution function $\Frm$ in \eqref{eq:pv.model-based.2} equals $\Frm_{p,n-p}$, so $\pi(\y)$ coincides with \eqref{eq:p-value.model.F} in the special case of $p_o = 0$. This follows from a standard argument for linear models: Let $\bb_1,\bb_2,\ldots,\bb_n$ be an orthonormal basis of $\R^n$ such that $\mathrm{span}(\x_1,\ldots,\x_p) = \mathrm{span}(\bb_1,\ldots,\bb_p)$. Then $\z \sim N_n(\bs{0},\bs{I})$ has the same distribution as $\tilde{\z} = \sum_{i=1}^n z_i \bb_i$, and
\[
	S(\tilde{\z}) \ = \ \frac{\sum_{j=1}^p z_j^2 / p}{\sum_{i=p+1}^n z_i^2 / (n - p)}
\]
has distribution function $\Frm_{p,n-p}$ by definition of Fisher's F distributions.
\end{exmp}

\begin{exmp}[Multiple t tests]
\label{ex:multiple.t}
Suppose that the linear span $\bbV$ of the regressors $\x_1,\ldots,\x_p$ satisfies $q := \dim(\bbV) < n$. Further let $\mathbb{A}$ be a subset of $\bbV \cap \bbS_n$. With the orthogonal projection $\hat{\y}$ of $\y$ onto $\bbV$, a possible test statistic is given by
\begin{equation}
\label{eq:T.statistic.1}
	S(\y) \ := \ \hat{\sigma}^{-1} \sup_{\bs{a} \in \mathbb{A}} \, |\bs{a}^\top \y|
	\quad\text{with}\quad
	\hat{\sigma} \ := \ (n - q)^{-1/2} \|\y - \hat{\y}\| .
\end{equation}
Note that under $H_o$, each term $\hat{\sigma}^{-1} \bs{a}^\top\y$ follows student's t distribution with $n-q$ degrees of freedom. This example of $S(\cdot)$ is motivated by Tukey's studentized maximum modulus or studentized range test statistics; see \cite{Miller_1981}.
\end{exmp}

\begin{exmp}[Multiple F tests]
\label{ex:multiple.F}
Let $p$ and $\x_1,\ldots,\x_p$ be arbitrary, and let $\Lambda$ be a family of subsets $M$ of $\{1,\ldots,p\}$ such that the vectors $\x_j$, $j \in M$, are linearly independent with $\# M < n$. With $\Pi_M$ denoting the orthogonal projection from $\R^n$ onto $\mathrm{span}(\x_j : j \in M)$, a possible test statistic is given by
\begin{equation}
\label{eq:T.statistic.2}
	S(\y) \ := \ \max_{M \in \Lambda} \,
		\frac{\|\Pi_M\y\|^2/\# M}{\|\y - \Pi_M\y\|^2/(n - \# M)} .
\end{equation}
The idea behind this test statistic is that possibly $\y = \bmu + \beps$ with a random vector $\beps$ having spherically symmetric distribution and a fixed vector $\bmu \in \R^n$ such that
\[
	\|\Pi_M\bmu\|^2
	\ \gg \ \| \bmu - \Pi_M\bmu\|^2
\]
for some $M \in \Lambda$.
\end{exmp}

%----------------------------------------
\subsection{The model-free point of view}
%----------------------------------------

To elaborate on the connection between model-based and model-free approach, note first that the null hypothesis $H_o$ is equivalent to an orthogonal invariance property. With $\stackrel{d}{=}$ denoting equality in distribution, the alternative formulation reads as follows.

\noindent
$H_o'$: \ $\T\y \stackrel{d}{=} \y$ for any fixed $\T \in \bbO_n$.

Another equivalent formulation involves normalized Haar measure $\Haar_n$ on $\bbO_n$. This is the unique distribution of a random matrix $\bH \in \bbO_n$ with left-invariant distribution in the sense that
\[
	\T\bH \ \stackrel{d}{=} \ \bH \quad\text{for any fixed} \ \T \in \bbO_n .
\]
For a thorough account of Haar measure we refer to \cite{Eaton_1989}; in Appendix~\ref{app:Haar} we mention two explicit constructions of $\bH$ and resulting properties. For the moment it suffices to know that also
\[
	\bH^\top \ \stackrel{d}{=} \ \bH \ \stackrel{d}{=} \ \bH\T 
	\quad\text{for any fixed} \ \T \in \bbO_n .
\]
Moreover, for any fixed unit vector $\bv \in \bbS_n$, the random vector $\bH\bv$ is uniformly distributed on $\bbS_n$. Now the null hypothesis $H_o'$ may be reformulated as follows:

\noindent
$H_o''$: \ If $\bH \sim \Haar_n$ is independent from $\y$, then $\bH\y \stackrel{d}{=} \y$.

\noindent
The equivalence of the null hypotheses $H_o$, $H_o'$ and $H_o''$ is explained in Appendix~\ref{app:Equivalence.hypotheses}.

From now on suppose that the test statistic $S(\y,\x_1,\ldots,\x_p)$ is orthogonally invariant in the sense that
\begin{equation}
\label{eq:ortho-invariance}
	S(\T\y, \T\x_1, \ldots, \T\x_p) \ = \ S(\y, \x_1, \ldots, \x_p)
	\quad\text{for all} \ \y \in \R^n \ \text{and} \ \T \in \bbO_n .
\end{equation}
Since $\x \mapsto \T\x$ preserves inner products, a sufficient condition for orthogonal invariance of the test statistic $S$ is that $S(\y,\x_1,\ldots,\x_p)$ depends only on the inner products $\y^\top \y$, $\y^\top \x_j$ and $\x_j^\top \x_k^{}$, $1 \le j \le k \le p$. Then the p-value \eqref{eq:pv.model-based.1} may be rewritten as follows:
\begin{align*}
	\pi(\y) \
	&= \ \Pr \bigl( S(\|\y\|\bu, \x_1,\ldots,\x_p) \ge S(\y) \,\big|\, \y \bigr) \\
	&= \ \Pr \bigl( S(\bH\y, \x_1,\ldots,\x_p) \ge S(\y) \,\big|\, \y \bigr) \\
	&= \ \Pr \bigl( S(\y, \bH^\top\x_1,\ldots,\bH^\top\x_p) \ge S(\y) \,\big|\, \y \bigr) \\
	&= \ \Pr \bigl( S(\y, \bH\x_1,\ldots,\bH\x_p) \ge S(\y,\x_1,\ldots,\x_p)
		\,\big|\, \y \bigr) ,
\end{align*}
where $\bH \sim \Haar_n$ is independent from $\y$. If we adopt the model-free point of view and consider all vectors $\y$ and $\x_1,\ldots,\x_p$ as fixed, we may write
\[
	\pi(\y)
	\ = \ \Pr \bigl( S(\y, \bH\x_1,\ldots,\bH\x_p) \ge S(\y,\x_1,\ldots,\x_p) \bigr) .
\]
Thus, $\pi(\y)$ measures the strength of the apparent association between $\y$ and the regressor tuple $(\x_1,\ldots,\x_p)$, as quantified by the test statistic $S(\y,\x_1,\ldots,\x_p)$, by comparing the latter value with $S(\y,\bH\x_1,\ldots,\bH\x_p)$. That means, the regressor tuple $(\x_1,\ldots,\x_p)$ undergoes a random orthogonal transformation, and there is certainly no ``true association'' between $\y$ and $(\bH\x_1,\ldots,\bH\x_p)$. To make the latter point rigorous, note that if $\bH, \bs{J} \sim \Haar_n$ and $\bu \sim \mathrm{Unif}(\bbS_n)$ are independent (while $\y$ is fixed), then $\bs{J}^\top\bH \sim \Haar_n$ too, whence
\begin{align*}
	S(\y, \bH\x_1, \ldots, \bH\x_p) \
	&\stackrel{d}{=} \
		S(\y, \bs{J}^\top\bH\x_1, \ldots, \bs{J}^\top\bH\x_p) \\
	&= \ S(\bs{J}\y, \bH\x_1, \ldots, \bH\x_p) \\
	&\stackrel{d}{=} \ S(\|\y\| \bu, \bH\x_1, \ldots, \bH\x_p) .
\end{align*}

Finally, recall that the method in the introduction with $p_o = 0$ amounts to replacing the fixed regressors $\x_1,\ldots,\x_p$ with independent random vectors $\x_1^*,\ldots,\x_p^* \sim N_n(\bs{0},\bs{I})$. But in connection with the F test statistic, this has the same effect as replacing the former random vectors with $\bH\x_1,\ldots,\bH\x_p$. Indeed, in case of linearly independent vectors $\x_1,\ldots,\x_p$, the value of $S(\y)$ in \eqref{eq:F.statistic} depends only on $\y$ and the $p$-dimensional linear space $\mathrm{span}(\x_1,\ldots,\x_p)$. Moreover, the distributions of $\mathrm{span}(\bH\x_1,\ldots,\bH\x_p)$ and of $\mathrm{span}(\x_1^*,\ldots,\x_p^*)$ coincide, see Appendix~\ref{app:Haar}.

%-----------------------------------
\subsection{Composite null models}
\label{subsec:Composite.null.models}
%-----------------------------------

Quite often, the potential influence of some regressors $\x_1,\ldots,\x_{p_o}$ with $1 \le p_o < p$ is out of question or not of primary interest, and the main question is whether the regressors $\x_{p_o+1},\ldots,\x_p$ are really relevant for the approximation of $\y$. One possibility to deal with that is to ``residualize'' the response $\y$ and the regressors $\x_{p_o+1}, \ldots, \x_p$, that is, to project them onto the orthogonal complement of $\mathrm{span}(\x_1,\ldots,\x_{p_o})$. In case of $p_o = 1$ and $\bs{x}_1 = (1)_{i=1}^n$, this boils down to centering $\y$ and $\bs{x}_2,\ldots,\bs{x}_p$.

More formally, assuming without loss of generality that $\x_1,\ldots,\x_{p_o}$ are linearly independent, let $\bb_1,\bb_2,\ldots,\bb_n$ be an orthonormal basis of $\R^n$ such that $\bb_1,\ldots,\bb_{p_o}$ form a basis of $\spann(\x_1,\ldots,\x_{p_o})$. With $\bB = [\bb_{p_o+1},\ldots,\bb_n] \in \R^{n\times (n-p_o)}$, the model equation \eqref{eq:Gaussian.model} implies that
\[
	\bB^\top\y \
	\sim \ N_{n - p_o} \Bigl( \sum_{j=p_o+1}^p \beta_j \bB^\top\x_j , \sigma^2 \bs{I} \Bigr) .
\]
Generally, the previous model-based and model-free considerations can be applied to $\bB^\top \y \in \R^{n-p_o}$ in place of $\y$.

Applying the model-based or model-free approach with the F test statistic \eqref{eq:F.statistic} to the transformed observations $\bB^\top\y$ and $\bB^\top\x_{p_o+1}, \ldots, \bB^\top \x_p$ yields the p-value \eqref{eq:p-value.model.beta} in the introduction, see also the first part of the proof of Lemma~\ref{lem:p-value.model-free.beta}.

%------------------------------
\section{Equivalence regions}
\label{sec:Equivalence.regions}
%------------------------------

%----------------------------------
\subsection{General considerations}
%----------------------------------

\paragraph{Model-based approach.}
Let $\bbM \subset \R^n$ be a given set. We assume that $\y$ is a random vector such that
\begin{equation}
\label{eq:standard.model.CR}
	\y \ \stackrel{d}{=} \ \bmu + \beps
\end{equation}
with an unknown fixed parameter vector $\bmu \in \bbM$ and a random vector $\beps$ with spherically invariant distribution on $\R^n$, where $\Pr(\beps = \bs{0}) = 0$. Now let $S(\y) = S(\y,\x_1,\ldots,\x_p)$ be a test statistic which is scale-invariant in $\y \ne \bs{0}$, and let $S(\bs{0}) := 0$. For a given (small) number $\alpha \in (0,1)$ we define the equivalence region
\[
	C_\alpha(\y) = C_\alpha(\y,\x_1,\ldots,\x_p)
	\ := \ \bigl\{ \bs{m} \in \bbM : S(\y - \bs{m}) \le \kappa_\alpha \bigr\} ,
\]
where $\kappa_\alpha$ is the $(1 - \alpha)$-quantile of the distribution of $S(\bu) \stackrel{d}{=} S(\z)$ with random vectors $\bu \sim \mathrm{Unif}(\bbS_n)$ and $\z \sim N_n(\bs{0},\bs{I})$. This defines an $(1 - \alpha)$-confidence region for $\bmu$ in the sense that in case of \eqref{eq:standard.model.CR},
\[
	\Pr(C_\alpha(\y) \ni \bmu)
	\ = \ \Pr(S(\z) \le \kappa_\alpha)
	\ \ge \ 1 - \alpha .
\]

\paragraph{Model-free interpretation.}
We consider $\y$ as fixed and assume in addition that the test statistic $S(\y,\x_1,\ldots,\x_p)$ is orthogonally invariant. Then the equivalence region $C_\alpha(\y)$ consists of all vectors $\bs{m} \in \bbM$ such that the association between $\y-\bs{m}$ and the tuple $(\x_1,\ldots, \x_p)$ is not substantially stronger than the association between $\y-\bs{m}$ and the randomly rotated tuple $(\bH\x_1,\ldots,\bH\x_p)$, where $\bH \sim \Haar_n$. Precisely, the value $S(\y - \bs{m},\x_1,\ldots,\x_p)$, our measure of association, is not larger than the $(1 - \alpha)$-quantile of $S(\y - \bs{m},\bH\x_1,\ldots,\bH\x_p)$.

\paragraph{Example~\ref{ex:F} (continued)}
Let $\bbM$ be the linear space $\mathrm{span}(\x_1,\ldots,\x_p) = \X\R^p$ with the matrix $\X = [\x_1,\ldots,\x_p] \in \R^{n\times p}$. Then the equivalence region equals
\[
	C_\alpha(\y) \ = \ \bigl\{ \X\bs{\eta} : \bs{\eta} \in \R^p ,
		\|\hat{\y} - \X\bs{\eta}\|^2 \le p F_{p,n-p}^{-1}(1 - \alpha) \hat{\sigma}^2(\y)
		\bigr\} ,
\]
where $\hat{\sigma}^2(\y) = \|\y - \hat{\y}\|^2/(n-p)$. The corresponding set $\tilde{C}_\alpha(\y) = \bigl\{ \bs{\eta} \in \R^p : \X\bs{\eta} \in C_\alpha(\y) \bigr\}$ is Scheff\'{e}'s well-known confidence ellipsoid for the unknown parameter $\bs{\beta} \in \R^p$ such that $\bmu = \X\bs{\beta}$.

%----------------------------------------
\subsection{Inference on a sparse signal}
\label{subsec:Sparsity}
%----------------------------------------

Suppose that $p = n$, and that the vectors $\x_1,\ldots,\x_n$ are linearly independent. In that case, $\X := [\x_1,\ldots,\x_n]$ is nonsingular, and assuming \eqref{eq:standard.model.CR}, the least squares estimator of $\bs{\beta} := \X^{-1}\bmu$ is given by $\hat{\bs{\beta}}(\y) := \X^{-1}\y$. Writing $\X^{-1} = [\ba_1,\ldots,\ba_n]^\top$, the Gauss--Markov estimator of $\beta_i$ equals $\hat{\beta}_i(\y) = \ba_i^\top\y$. In case of $\y \sim N_n(\bmu,\sigma^2\bs{I})$, it has distribution $N(\beta_i, \|\ba_i\|^2 \sigma^2) = N \bigl( \beta_i, ((\X^\top\X)^{-1})_{ii}^{} \sigma^2 \bigr)$.

Suppose that $\bmu = \X\bs{\beta}$ is sparse in the sense that
\[
	\|\bs{\beta}\|_0 \ := \ \#\{i \le n : \beta_i \ne 0\}
\]
is relatively small compared with $n$. Then a possible test statistic for the null hypothesis ``$\bmu = \bs{0}$'' is given by
\[
	S_\ell(\y) = S_\ell(\y,\x_1,\ldots,\x_n)
	\ := \ G_{1}(\y) / G_{\ell}(\y)
\]
for some integer $\ell \in \{2,\ldots,n\}$, where $G_{1}(\y) \ge G_{2}(\y) \ge \cdots \ge G_{n}(\y)$ are the values $\|\ba_1\|^{-1}|\ba_1^\top\y|$, $\|\ba_2\|^{-1}|\ba_2^\top\y|$, \ldots, $\|\ba_n\|^{-1}|\ba_n^\top\y|$ in descending order. The idea behind this choice is that $G_{\ell}(\y)$ depends mostly on the noise vector $\beps$ if $\|\bs{\beta}\|_0 < \ell$, while $G_{1}(\y)$ is mainly driven by $\bmu$ in case of $\|\bs{\beta}\| \gg 0$. Scale-invariance of the test statistic $S_\ell(\y,\x_1,\ldots, \x_n)$ in $\y$ is obvious. Since $\X^{-1} = (\X^\top\X)_{}^{-1} \X^\top$, one may write $\|\ba_i\|^2 = ((\X^\top\X)^{-1})_{ii}$ and $\bs{a}_i^\top\y = \bigl( (\X^\top\X)^{-1} \X^\top \y \bigr)_i$, whence $S_\ell(\cdots)$ is also orthogonally invariant.

Denoting the $(1 - \alpha)$-quantile of $S_\ell(\z)$, $\z \sim N_n(\bs{0},\bs{I})$, with 
$\kappa_{\ell,\alpha}$ and setting $\bbM = \R^n$, we obtain the equivalence region
\[
	C_{\ell,\alpha}(\y)
	\ := \ \bigl\{ \bs{m} \in \R^n :
		S_\ell(\y - \bs{m}) \le \kappa_{\ell,\alpha} \bigr\} .
\]
This region is rather useless per se. But if we restrict our attention to vectors $\bs{m} = \X\bb$, where $\|\bb\|_0$ is bounded by a given number, we end up with the equivalence regions
\[
	\tilde{C}_{k,\ell,\alpha}(\y)
	\ := \ \bigl\{ \bb \in \R^n : \|\bb\|_0 \le k ,
		S_\ell(\y - \X\bb) \le \kappa_{\ell,\alpha} \bigr\} ,
		\quad 1 \le k < n ,
\]
which are potentially useful in case of $k < \ell$. In this manuscript we only prove a first result about these equivalence regions.

\begin{lem}
\label{lem:Sparse.ER}
Let
\[
	k_{\ell,\alpha}(\y) \ := \ \ell - 1 - \max \bigl\{ j - i : 1 \le i \le j \le \ell,
		G_i(\y)/G_j(\y) \le \kappa_{\ell,\alpha} \bigr\} .
\]
Then $\tilde{C}_{\alpha,k,\ell}(\y) \ne \emptyset$ if and only if $k \ge k_{\ell,\alpha}(\y)$.
\end{lem}

In the model-based context, $k_{\ell,\alpha}$ is a lower $(1 - \alpha)$-confidence bound for $\|\bs{\beta}\|_0$, and $\tilde{C}_{k,\ell,\alpha}$ is a $(1 - \alpha)$-confidence region for $\bs{\beta}$ under the additional assumption that $\|\bs{\beta}\|_0 \le k$.

%-----------------------------------------------------
\subsection{Special case: the Gaussian sequence model}
\label{subsec:Sequence}
%-----------------------------------------------------

Suppose that $\x_1,\ldots,\x_n$ is the standard basis of $\R^n$. Then $(G_i(\z))_{i=1}^n$ has the same distribution as $\bigl( \tilde{\Phi}^{-1}(U_{(n+1-i)}) \bigr)_{i=1}^n$, where $U_{(1)} < \cdots < U_{(n)}$ are the order statistics of independent random variables $U_1,\ldots,U_n \sim \mathrm{Unif}[0,1]$, and $\tilde{\Phi}$ is the distribution function of $|z_1|$, that is, $\tilde{\Phi}(r) = 2\Phi(r) - 1$ for $r \ge 0$. Since $U_{(n+1-\ell)} = 1 - \ell/n + O_p(n^{-1/2})$ as $n \to \infty$, a reasonable choice for $\ell$ seems to be $\ell \approx (1 - \tilde{\Phi}(1)) n = 2 \Phi(-1) n \approx 0.32 \, n$. Then $\kappa_{\ell,\alpha}$ is approximately the $(1 - \alpha)$-quantile of $G_1(\z) = \max_{1 \le i \le n} |z_i| = \sqrt{2 \log n} (1 + o_p(1))$.

Specifically, let $n = 100$ and $\ell = 32$. Numerical computations outlined in Appendix~\ref{app:Sequence} yield $\kappa_{\ell,0.01} = 4.083$.

Now we consider the model-based setting with $\y \sim N_n(\bmu,\bs{I})$ and and two different choices for $\bmu$. The proof of Lemma~\ref{lem:Sparse.ER} includes the construction of a particular vector $\hat{\bmu}(\y) = \hat{\bmu}_{\ell,\alpha}(\y)$ in $C_{\ell,\alpha}(\y)$ such that $\|\hat{\bmu}\|_0 = k_{\ell,\alpha}(\y)$. Now we investigated the distribution of the latter number, of the set $\{i \le n : \hat{\mu}_i \ne 0\}$ and of the cosine of the angle between $\bmu$ and $\hat{\bmu}$.

In the first scenario, we considered the sparse vector $\bmu = (10,-6,3,0,\ldots,0)^\top$. In 100'000 Monte Carlo simulations we estimated the joint distribution of
\[
	\mathrm{TP}(\y) \ := \ \bigl\{ i \in \{1,2,3\} : \hat{\mu}_i(\y) \ne 0 \bigr\}
	\quad\text{and}\quad
	\mathrm{FP}(\y) \ := \ \# \{i > 3 : \hat{\mu}_i(\y) \ne 0 \} ,	
\]
the \textsl{set} of nonzero components of $\bmu$ which were detected correctly (``true positives'') and the \textsl{number} of ``false positives'', respectively. Table~\ref{tab:MCa} contains estimated probabilities rounded to four digits. It turned out that $\mathrm{FP}(\y) \ge 2$ with estimated probability less than $10^{-4}/2$, and $\mathrm{FP}(\y) \ge 1$ with estimated probability $0.0058$. The first component of $\bmu$ was always identified correctly as non-zero, whereas the second and third component stayed sometimes undetected.

\begin{figure}
\centering
\includegraphics[width=0.75\textwidth]{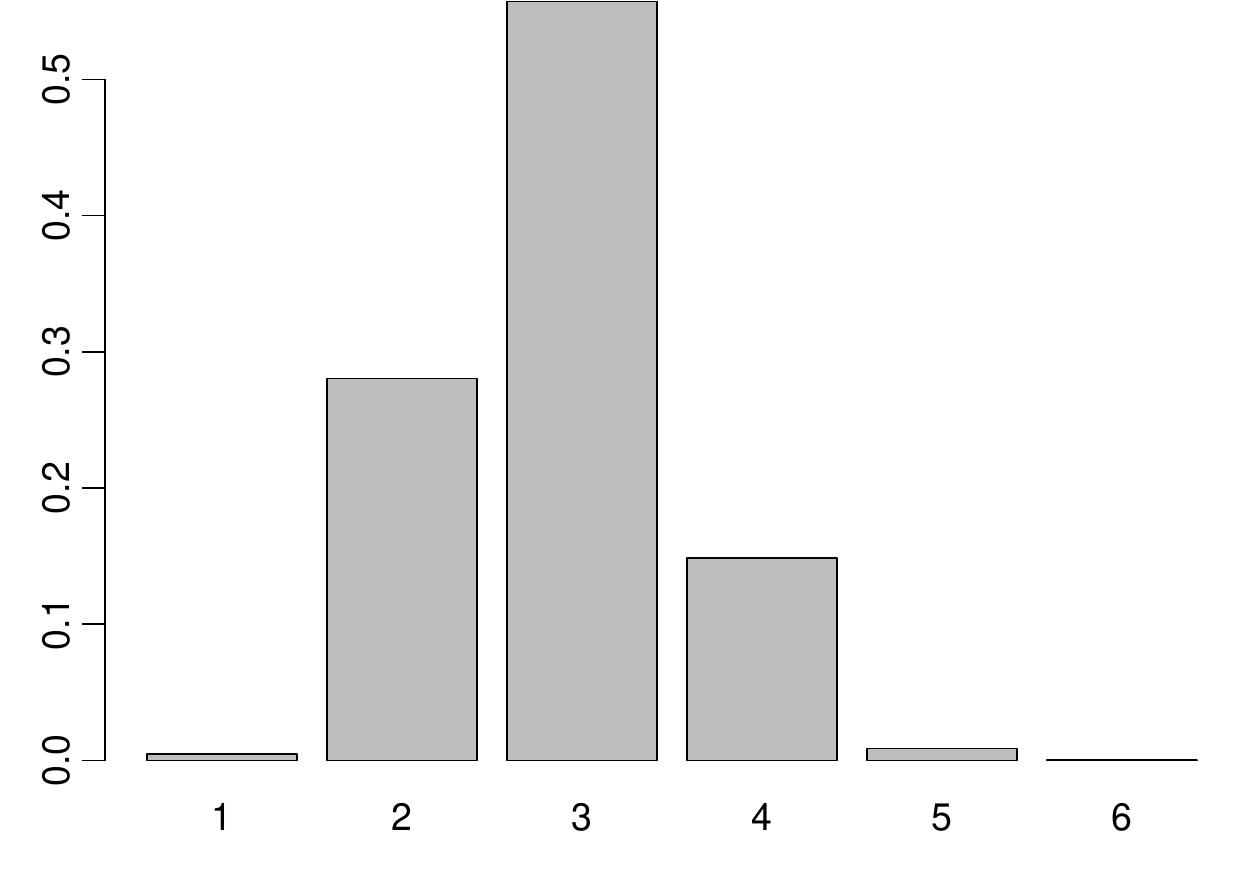}
\caption{Distribution of $k_{\ell,\alpha}(\y)$ in 2nd scenario.}
\label{fig:MCb_k}
\end{figure}

In the second scenario, we considered the non-sparse vector $\bmu = \bigl( (-1)^{i-1} 2^{-(i-1)/2} \bigr)_{i=1}^n$. Although $\|\bmu\|_0 = n$, the first four to six components of $\bmu$ contain the main signal, because $\|\bmu\|^{-2} \sum_{i > 4} \mu_i^2 = 0.0625$ and $\|\bmu\|^{-2} \sum_{i > 6} \mu_i^2 < 0.0157$. Figure~\ref{fig:MCb_k} shows the (estimated) distribution of $k_{\ell,\alpha}(\y)$. Recall that the latter is a lower $(1 - \alpha)$-confidence bound for $\|\bs{\beta}\|_0$.

\begin{figure}
\includegraphics[width=0.45\textwidth]{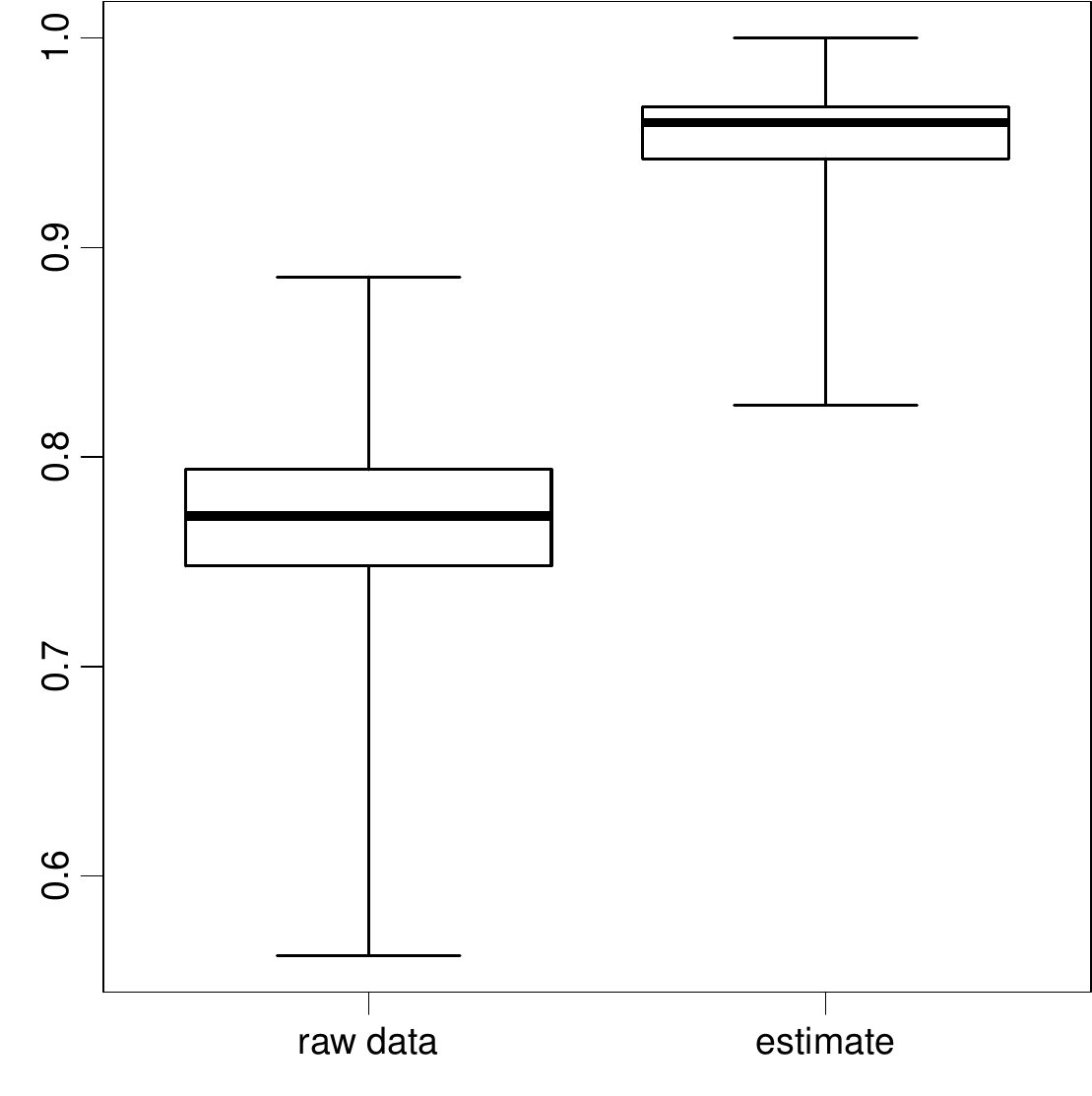}
\hfill
\includegraphics[width=0.45\textwidth]{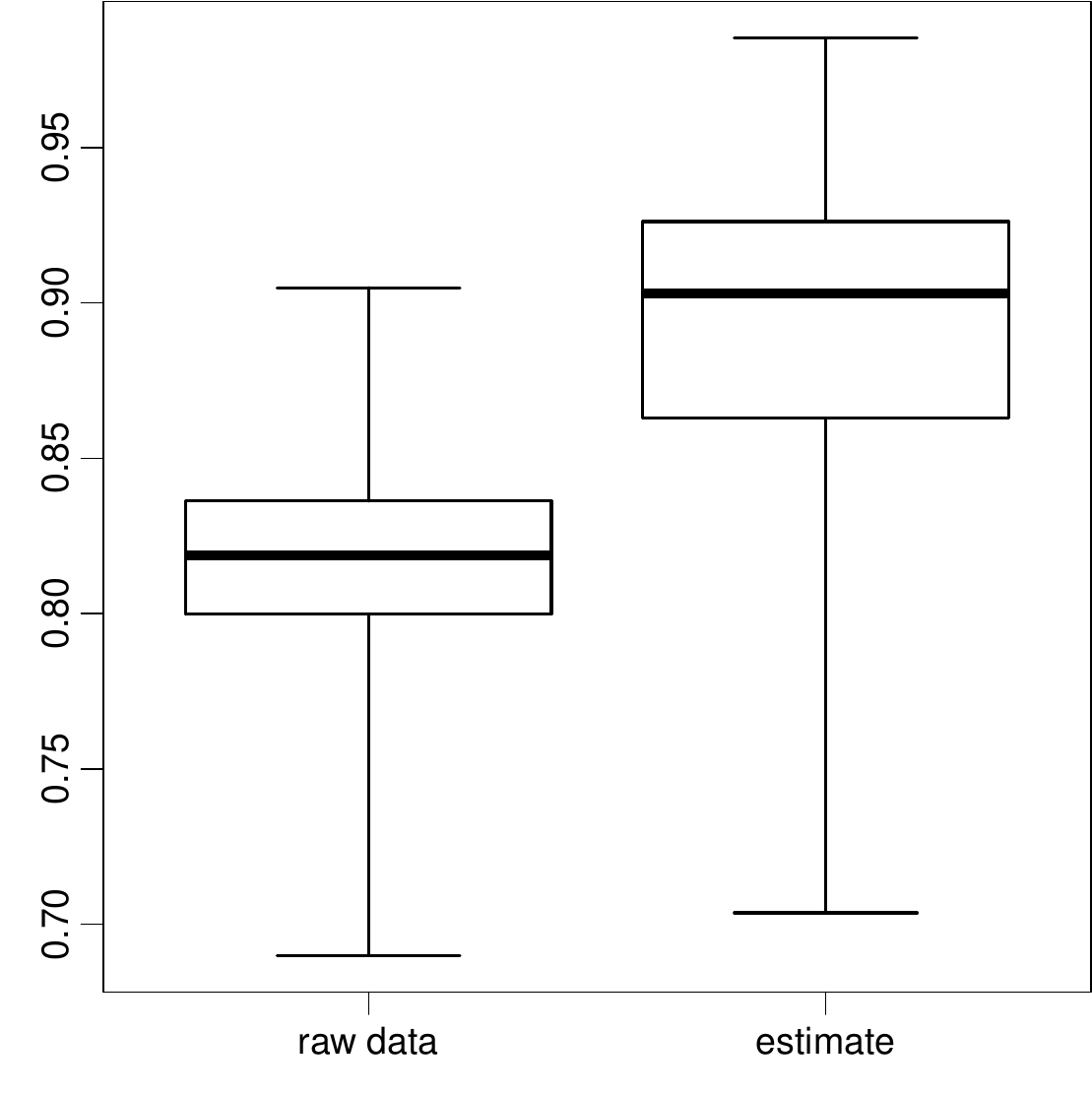}
\caption{Distribution of $\cos \angle(\bmu,\y)$ and $\cos \angle(\bmu, \hat{\bmu})$ in 1st scenario (left) and 2nd scenario (right).}
\label{fig:MC_cos}
\end{figure}

Finally, Figure~\ref{fig:MC_cos} show boxplots of the distribution of
\[
	\cos \angle(\bmu,\y)
	\quad\text{and}\quad
	\cos \angle(\bmu, \hat{\bmu})
\]
for both scenarios. These plots illustrate that the sparse estimator $\hat{\bmu}$ captures $\bmu$ better than the raw data $\y$.

\begin{table}
\caption{Joint distribution of $\mathrm{TP}(\y)$ and $\mathrm{FP}(\y)$ in 1st scenario (in per cent).}
\label{tab:MCa}
\[
	\begin{array}{lcccccc}
	\hline
		\multicolumn{1}{l}{\mathrm{TP}(\y) =}
		& \{1,2,3\} & \{1,2\} & \{1,3\} &
		\{1\} & \text{otherwise} & \text{total} \\
		% &  \{2,3\}, \{2\}, \{3\}, \emptyset & \\
	\hline\hline
	\mathrm{FP}(\y) = 0
		& 11.13 & 82.70 & 0.35 & 5.24 & 0.00 & 99.42 \\
	\hline
	\mathrm{FP}(\y) = 1
		&  0.12 &  0.42 & 0.00 & 0.04 & 0.00 &  0.58 \\
	\hline\hline
	\text{total}
		& 11.25 & 83.12 & 0.35 & 5.28 & 0.00 & 100.00 \\
	\hline
	\end{array}
\]
\end{table}

%----------------------------------------
\section{Final comments and further work}
\label{sec:Final.comments}
%----------------------------------------

\paragraph{Permutation tests.}
Let us relate the considerations in the present paper to research about permutation tests in general regression models. We start from the model-based approach with $\y = \bs{X} \bs{\beta} + \beps$ with fixed $\bs{X} = [\x_1,\ldots,\x_p] \in \R^{n\times p}$, $\bs{\beta} \in \R^p$ and a random error $\beps \in \R^n$. When testing the null hypothesis that $(\beta_j)_{p_o < j \le p} = \bs{0}$, the goal is to relax the assumption of a spherically symmetric distribution of $\beps$ to assuming exchangeability only. That is, for any fixed permutation $\tau$ of $\{1,\ldots,n\}$, the distribution of $\tau(\beps)$ coincides with the distribution of $\beps$, where $\tau(\bs{v}) := (v_{\tau(i)})_{i=1}^n$ for $\bs{v} \in \R^n$. Indeed, if $p_o = 0$ or $p_o = 1$ and $\bs{x}_1 = (1)_{i=1}^n$, the null hypothesis can be tested exactly with a permutation test in which the original data $(\y,\bs{X})$ are compared with $(\tau(\y), \bs{X})$ for all (or many randomly chosen) permutations $\tau$. In other cases, however, it is not obvious how to perform a valid permutation test. \cite{Kennedy_1995} describes potential pitfalls, and \cite{Winkler_etal_2014} provide a good overview of proposed permutation schemes. Apart from the aforementioned simple setting, all these proposals are justified by asymptotic considerations only. A particularly interesting reference is \cite{Freedman_Lane_1983}, because they motivated their permutation test also by the wish to interpret the resulting p-values in a model-free way. Their approach works as follows: Let $\hat{\beps}_o := \y - \hat{\y}_o$ be the residual vector under the null model. Then they propose to compare the F test statistic for the original observations $(\y,\bs{X})$ with the F test statistic applied to $(\hat{\y}_o + \tau(\hat{\beps}_o), \bs{X})$. The collection of data $(\hat{\y}_o + \tau(\hat{\beps}_o), \bs{X})$, where $\tau$ is an arbitrary permutation of $\{1,\ldots,n\}$, serves as a data-generated reference set of data to be compared with the original $(\y,\bs{X})$. If the latter sticks out in terms of the F test statistic, this is viewed as model-free evidence that the regressors $\x_j$, $p_o < j \le p$, are relevant. \cite{Freedman_Lane_1983} show that under certain regularity assumptions, the resulting p-value is asymptotically equivalent to the classical p-value \eqref{eq:p-value.model.F}. Thus our Lemma~\ref{lem:p-value.model-free.beta} may be viewed as a computationally simple and exact alternative to the paradigm of \cite{Freedman_Lane_1983}.

\paragraph{Variable selection.}
Sections~\ref{subsec:Sparsity} and \ref{subsec:Sequence} describe an ad hoc approach to model-free variable selection. Unfortunately it is computationally intensive, at least for non-orthogonal regressors $\x_1,\ldots,\x_n$. A much simpler and computationally very efficient method for variable selection has been developed and implemented by \cite{Davies_2024} as an R package \citep{R21}. A key element is a forward selection procedure, and we describe only the very first step in which one has to decide whether at least one of the covariate vectors $\x_1,\ldots,\x_p$ should be used to approximate $\y$. A possible test statistic for this task is given by
\[
	S(\y, \x_1,\ldots,\x_p) \ := \ \max_{j=1,\ldots,p} \frac{(\x_j^\top \y)^2}{\|\x_j\|^2 \|\y\|^2} ,
\]
which is just a monotone transformation of the multiple F test statistic in Example~\ref{ex:multiple.F} in the special case of $\Lambda = \bigl\{ \{1\}, \ldots, \{p\} \bigr\}$. There exists no closed formula for the model-free p-value
\[
	\Pr \bigl( S(\y, \bH\x_1,\ldots,\bH\x_p) \ge S(\y, \x_1,\ldots,\x_p) \bigr) ,
\]
so \cite{Davies_2024} uses the alternative p-value
\begin{align*}
	\Pr \bigl( & S(\y, \x_1^*,\ldots,\x_p^*) \ge S(\y, \x_1,\ldots,\x_p) \bigr) \\
	&= \ 1 - \Pr \biggl( \frac{(\x_j^{*\top}\y)^2}{\|\x_j^*\|^2 \|\y\|^2} < S(\y, \x_1,\ldots,\x_p) \
		\text{for} \ 1 \le j \le p \biggr) \\
	&= \ 1 - \mathrm{Beta}_{1/2,(n-1)/2} \bigl( S(\y, \x_1,\ldots,\x_p) \bigr)^p .
\end{align*}
If this modified p-value is smaller than a given threshold $\alpha$, one chooses a covariate vector $\x_{j_o}$ maximizing $\|\x_j\|^{-2} \|\y\|^{-2} (\x_j^\top\y)^2$ over all $j \in \{1,\ldots,p\}$ to approximate $\y$. Thereafter one replaces $\y$ and the remaining covariate vectors $\x_j$, $j \ne j_o$, with their orthogonal projections onto $\x_{j_o}^\perp$ and applies an analogous procedure\ldots

The price to be paid for the simplified p-value is that potential collinearities between the covariate vectors $\x_1,\ldots,\x_p$ are not taken into account, so the selection procedure is possibly conservative in the sense of selecting too few variables. But this method of ``Gaussian covariates'' yields satisfactory results even in very high-dimensional settings where alternative proposals fail. It can also be applied to vector autoregression.

\paragraph{Acknowledgements.}
The authors are grateful to a reviewer, Christian Hennig (editor) and Jan Hannig for constructive comments. Part of this work was supported by Swiss National Science Foundation.

%--------
\appendix

%---------------------------------------------------------
\section{Proof of Lemma~\ref{lem:p-value.model-free.beta}}
\label{app:proof.F.lemma}
%---------------------------------------------------------

Note that $\y$ is the sum of the three orthogonal vectors $\hat{\y}_o$, $\hat{\y}^* - \hat{\y}_o$ and $\y - \hat{\y}^*$. Thus,
\[
	\frac{\|\y - \hat{\y}^*\|^2}{\|\y - \hat{\y}_o\|^2} \
	= \ \frac{\|\y - \hat{\y}_o\|^2 - \|\y - \hat{\y}^*\|^2}{\|\y - \hat{\y}_o\|^2} \
	= \ 1 - \frac{\|\hat{\y}^* - \hat{\y}_o\|^2}{\|\y - \hat{\y}_o\|^2} ,
\]
whence the claim is equivalent to
\[
	\frac{\|\hat{\y}^* - \hat{\y}_o\|^2}{\|\y - \hat{\y}^*\|^2} \
	\sim \ \mathrm{Beta} \bigl( (p-p_o)/2, (n - p)/2 \bigr) .
\]

Let us first recall two well-known facts about a standard Gaussian random vector $\z = (z_i)_{i=1}^n$ in $\R^n$:

\noindent(F1) \ For any fixed matrix $\bB \in \mathbb{O}_n$, the random vector $\bB\z$ is standard Gaussian too. Equivalently, for any orthonormal basis $\bb_1,\ldots,\bb_n$ of $\R^n$, the random vector $\sum_{i=1}^n z_i \bb_i$ is standard Gaussian.

\noindent(F2) \ For any $k \in \{1,\ldots,n-1\}$, the random variable $\sum_{j=1}^k z_j^2 \big/ \sum_{i=1}^n z_i^2$ follows the beta distribution with parameters $k/2$ and $(n - k)/2$.

\paragraph{Step~1: Reduction to the case of $p_o = 0$.}
Suppose that $p_o \ge 1$. Let $\Pi_o$ be the orthogonal projection from $\R^n$ onto $\spann(\x_1,\ldots,\x_{p_o})$, so $\hat{\y}_o = \Pi_o\y$. The vector $\hat{\y}^* - \hat{\y}_o$ is the orthogonal projection of $\y - \Pi_o\y$ onto the linear span of the vectors $\x_j^* - \Pi_o\x_j^*$, $p_o < j \le p$. All these vectors lie in the $(n - p_o)$-dimensional linear space $\{\x_1,\ldots,\x_{p_o}\}^\perp$, and the independent random vectors $\x_j^* - \Pi_o\x_j^*$, $p_o < j \le p$, follow a standard Gaussian distribution on that space. The latter claim follows from (F1) when choosing an orthonormal basis of $\R^n$ such that $n - p_o$ of these basis vectors span $\{\x_1,\ldots,\x_{p_o}\}^\perp$. Hence, we may assume without loss of generality that $p_o = 0$ and $\hat{\y}_o = \bs{0}$.

\paragraph{Step~2: The case of $p_o = 0$.}
We have to show that $\|\y\|^{-2} \|\hat{\y}^*\|^2$ follows a beta distribution with parameters $p/2$ and $(n-p)/2$. To this end we argue similarly as \cite{Frankl_Maehara_1990}. With the random matrix $\bs{X} := [\x_1^*,\ldots,\x_p^*] \in \R^{n\times p}$, the orthogonal projection $\hat{\y}^*$ is given by $\hat{\y}^* = \bs{X} (\bs{\X}^\top \bs{X})^{-1} \bs{X}^\top \y$, and with the unit vector $\bs{v} := \|\y\|^{-1} \y$,
\[
	\|\y\|^{-2} \|\hat{\y}^*\|^2 \
	= \ \bs{v}^\top \bs{X} (\bs{X}^\top\bs{X})^{-1} \bs{X}^\top \bs{v} .
\]
The distribution of this random variable does not depend on $\bs{v}$. Indeed, for any other unit vector $\bs{w} \in \R^n$, let $\bs{B} \in \bbO_n$ such that $\bs{B} \bs{v} = \bs{w}$. Since $\tilde{\bs{X}} := \bs{B} \bs{X}$ has the same distribution as $\bs{X}$, and since $\bs{X}^\top\bs{X} = \tilde{\bs{X}}^\top \tilde{\bs{X}}$, we may conclude that
\[
	\bs{v}^\top \bs{X} (\bs{X}^\top\bs{X})^{-1} \bs{X}^\top \bs{v} \
	= \ \bs{w}^\top \tilde{\bs{X}} (\tilde{\bs{X}}^\top\tilde{\bs{X}})^{-1} \tilde{\bs{X}}^\top \bs{w} \
	\stackrel{d}{=} \ \bs{w}^\top \bs{X} (\bs{X}^\top \bs{X})^{-1} \bs{X}^\top \bs{w} .
\]
But then we may replace $\bs{v}$ with $\|\z\|^{-1} \z$ with a standard Gaussian random vector $\z$ which is independent from $\bs{X}$. Conditional on $\X$, this vector $\z$ has the same distribution as $\sum_{i=1}^n z_i \bs{b}_i$ with an orthonormal basis $\bb_1,\ldots,\bb_n$ of $\R^n$ such that $\spann(\x_1^*,\ldots,\x_p^*) = \spann(\bb_1,\ldots,\bb_p)$, and the orthogonal projection of $\sum_{i=1}^n z_i \bs{b}_i$ onto the latter space equals $\sum_{j=1}^p z_j \bb_j$. Consequently, conditional on $\bs{X}$,
\[
	\|\y\|^{-2} \|\hat{\y}^*\|^2
	\stackrel{d}{=} \ \Bigl\| \sum_{i=1}^n z_i \bb_i \Bigr\|^{-2}
		\Bigl\| \sum_{j=1}^p z_j \bb_j \Bigr\|^2 \
	= \ \sum_{j=1}^p z_j^2 \big/ \sum_{i=1}^n z_i^2
\]
follows a beta distribution with parameters $p/2$ and $(n-p)/2$. Since this does not depend on $\bs{X}$, we may conclude that the ratio $\|\y\|^{-2} \|\hat{\y}^*\|^2$ has distribution $\mathrm{Beta}(p/2, (n-p)/2)$.
\qed

%-----------------------------------------------
\section{Haar measure on $\bbO_n$ in a nutshell}
\label{app:Haar}
%-----------------------------------------------

Recall that $\Haar_n$ is the unique probability measure on $\bbO_n$ such that a random matrix $\bH \sim \Haar_n$ satisfies $\T\bH \stackrel{d}{=} \bH$ for any fixed $\T \in \bbO_n$ (left-invariance). For the reader's convenience, we explain a few well-known basic facts here.

\paragraph{Existence.}
To construct such a random matrix explicitly, let $\Z \in \R^{n\times n}$ be a random matrix such that $\T\Z \stackrel{d}{=} \Z$ for any fixed $\T \in \bbO_n$, and $\mathrm{rank}(\Z) = n$ almost surely. This is true, for instance, if the $n^2$ entries of $\Z$ are independent with distribution $N(0,1)$. Then one can easily verify that $\bH := \Z (\Z^\top \Z)^{-1/2}$ is a random orthogonal matrix with left-invariant distribution. Another possible construction would be to apply Gram-Schmidt orthogonalization to the columns $\z_1, \z_2, \ldots, \z_n$ of $\Z$. This leads to a random matrix $\bH \in \bbO_n$ with left-invariant distribution, because the coefficients for the Gram--Schmidt procedure involve only the inner products $\z_i^\top\z_j^{}$, and these remain unchanged if $\Z$ is replaced with $\T\Z$ for some $\T \in \bbO_n$.

\paragraph{Inversion-invariance and uniqueness.}
Let $\bs{J}, \bH$ be independent random matrices with left-invariant distribution on $\bbO_n$. By conditioning on $\bs{J}$ one can show that $\bs{J}^\top \bH \stackrel{d}{=} \bH$, and conditioning on $\bH$ reveals that $\bs{J}^\top \bH = (\bH^\top \bs{J})^\top \stackrel{d}{=} \bs{J}^\top$. Hence, $\bs{J}^\top \stackrel{d}{=} \bH$. In the special case that $\bs{J}$ is an independent copy of $\bH$, this shows that
\[
	\bH^\top \ \stackrel{d}{=} \ \bH .
\]
And then, with the original $\bs{J}$, we see that $\bs{J} = (\bs{J}^\top)^\top \stackrel{d}{=} \bH^\top \stackrel{d}{=} \bH$, i.e.
\[
	\bs{J} \ \stackrel{d}{=} \ \bH .
\]

\paragraph{Right-invariance.}
If $\bH \sim \Haar_n$, then its distribution is right-invariant in the sense that $\bH\T \stackrel{d}{=} \bH$ for any fixed $\T \in \bbO_n$. This follows easily from left-invariance and inversion-invariance.

\paragraph{Random linear subspaces of $\R^n$.}
The second construction of $\bH = [\bs{h}_1,\ldots,\bs{h}_n]$ via Gram--Schmidt may be applied to a random matrix $\Z$ with independent columns $\z_1, \z_2, \ldots, \z_n \sim N_n(\bs{0},\bs{I})$. This shows that the first column of $\bH$ has distribution $\mathrm{Unif}(\bbS_n)$. Combined with right-invariance, applied to permutation matrices, this shows that any column of $\bH$ has distribution $\mathrm{Unif}(\bbS_n)$. Furthermore, for $k \in \{1,2,\ldots,n\}$ and arbitrary linearly independent vectors $\x_1,\ldots,\x_k \in \R^n$,
\[
	\spann(\z_1,\ldots,\z_k)
	\ \stackrel{d}{=} \
	\spann(\bH_1,\ldots,\bH_k)
	\ \stackrel{d}{=} \
	\spann(\bH\x_1,\ldots,\bH\x_k) .
\]
The first equality follows from the construction of $\bH$. For the second one, let $\bB \in \R^{k\times k}$ be a nonsingular matrix such that the columns of the matrix $[\bs{t}_1,\ldots,\bs{t}_k] = [\x_1,\ldots,\x_k] \bB$ are orthonormal. Extending them to a orthonomal basis $\bs{t}_1,\bs{t}_2,\ldots,\bs{t}_n$ of $\R^n$,
\begin{align*}
	\spann(\bH\x_1,\ldots,\bH\x_k)
	&= \ \bigl\{ \bH [\x_1,\ldots,\x_k] \bs{\lambda} : \bs{\lambda} \in \R^k \bigr\} \\
	&= \ \bigl\{ \bH [\bs{t}_1,\ldots,\bs{t}_k] \bs{\lambda} : \bs{\lambda} \in \R^k \bigr\} ,
\end{align*}
and this is the linear span of the first $k$ columns of $\bH [\bs{t}_1,\ldots,\bs{t}_n]$. But by right-invariance, the latter matrix has the same distribution as $\bH$, because $[\bs{t}_1,\ldots,\bs{t}_n] \in \bbO_n$.

%---------------------------------------------------------------------
\section{Equivalence of the null hypotheses $H_o$, $H_o'$ and $H_o''$}
\label{app:Equivalence.hypotheses}
%---------------------------------------------------------------------

Recall that we consider a random vector $\y$ with continuous distribution on $\R^d$. Suppose first that $H_o$ is true. Then $\y$ has the same distribution as $\|\y\| \|\z\|^{-1} \z$, where $\y$ and $\z \sim N_n(\bs{0},\bs{I})$ are independent. For any fixed $\T \in \bbO_n$, orthogonal invariance of the standard Gaussian distribution implies that $\T\z \stackrel{d}{=} \z$, so
\[
	\T\y
	\ \stackrel{d}{=} \ \|\y\| \|\z\|^{-1} \T\z \ = \ \|\y\| \|\T\z\|^{-1} \T\z
	\ \stackrel{d}{=} \ \|\y\| \|\z\|^{-1} \z
	\ \stackrel{d}{=} \ \y .
\]
Consequently, $H_o'$ is true as well.

Now suppose that $H_o'$ is true. If $\bH \sim \Haar_n$ and $\y$ are independent, then for any measurable set $B \subset \R^n$,
\[
	\Pr(\bH\y \in B) \ = \ \Ex \Pr(\bH\y \in B \,|\, \bH)
	\ = \ \Ex \Pr(\y \in B \,|\, \bH) \ = \ \Pr(\y \in B) ,
\]
where $\Pr(\bH\y \in B \,|\, \bH) = \Pr(\Y \in B \,|\, \bH)$ by $H_o'$.

Finally, suppose that $H_o''$ is true. Then for any measurable set $B \subset \R^n$,
\begin{align*}
	\Pr(\y \in B)
	\ = \ \Pr(\bH\y \in B) \
	&= \ \Ex \Pr(\bH\y \in B \,|\, \y) \\
	&= \ \Ex \Pr \bigl( \|\y\| \bu \in B \,\big|\, \y \bigr)
	\ = \ \Pr \bigl( \|\y\| \bu \in B \bigr) ,
\end{align*}
where $\bu \sim \mathrm{Unif}(\bbS_n)$ and $\y$ are independent. Thus $H_o$ is satisfied as well.

%-------------------------------------------
\section{Proof of Lemma~\ref{lem:Sparse.ER}}
\label{app:proof.Lemma.sparsity}
%-------------------------------------------

We first construct a particular point $\hat{\bbeta} = \hat{\bbeta}_{\ell,\alpha}(\y) \in \R^n$ such that $S_\ell(\y - \X\hat{\bbeta}) \le \kappa_{\ell,\alpha}$ and $\|\hat{\bbeta}\|_0 = k_{\ell,\alpha}(\y)$. To this end let $(\sigma(1),\sigma(2),\ldots,\sigma(n))$ be a permutation of $(1,2,\ldots,n)$ such that $G_i(\y) = \|\ba_{\sigma(i)}\|^{-1} |\ba_{\sigma(i)}^\top\y|$ for $1 \le i \le n$. Let $k_{\ell,\alpha}(\y) = \ell - 1 - j_o + i_o$ with indices $1 \le i_o \le j_o \le \ell$ such that $|G_{i_o}(\y)|/|G_{j_o}(\y)| \le \kappa_{\ell,\alpha}$. Now we set
\[
	\hat{\beta}_{\sigma(i)} \ := \ \begin{cases}
		0
			&\text{if} \ i_o \le i \le j_o \ \text{or} \ i > \ell , \\
		\mathrm{sign}(\ba_{\sigma(i)}^\top\y)
			\bigl( |\ba_{\sigma(i)}^\top\y| - \|\ba_{\sigma(i)}\| G_{i_o}(\y) \bigr)
			&\text{if} \ 1 \le i < i_o , \\
		\mathrm{sign}(\ba_{\sigma(i)}^\top\y)
			\bigl( |\ba_{\sigma(i)}^\top\y| - \|\ba_{\sigma(i)}\| G_{j_o}(\y) \bigr)
			&\text{if} \ j_o < i \le \ell .
	\end{cases}
\]
Obviously, this defines a vector $\hat{\bbeta} \in \R^n$ such that $\|\hat{\bbeta}\|_0 = k_{\ell,\alpha}(\y)$. Note also that the numbers $\|\ba_{\sigma(i)}\|^{-1} \bigl|\ba_{\sigma(i)}^\top(\y - \X\hat{\bbeta}) \bigr|$ are nonincreasing in $i \in \{1,\ldots,n\}$ with
\[
	\|\ba_{\sigma(1)}\|^{-1} \bigl|\ba_{\sigma(1)}^\top(\y - \X\hat{\bbeta}) \bigr|
	\ = \ G_{i_o}(\y) , \quad
	\|\ba_{\sigma(\ell)}\|^{-1} \bigl|\ba_{\sigma(\ell)}^\top(\y - \X\hat{\bbeta}) \bigr|
	\ = \ G_{j_o}(\y) .
\]
Consequently, $\hat{\bbeta} \in \tilde{C}_{k,\ell,\alpha}(\y)$ for $k \ge k_{\ell,\alpha}(\y)$.

On the other hand, let $k_{\ell,\alpha}(\y) > 0$. We have to show that $S_\ell(\y - \bb) > \kappa_{\ell,\alpha}(\y)$ for any $\bb \in \R^n$ with $k := \|\bb\|_0 < k_{\ell,\alpha}(\y)$. Then there exist indices $1 \le i(1) < i(2) < \cdots < i(n-k) \le n$ and $1 \le j(1) < j(2) < \cdots < j(n-k) \le n$ such that $G_{i(s)}(\y) = G_{j(s)}(\y - \X\bb)$ for $1 \le s \le n-k$. But at least $n - k - (n - \ell) = \ell - k$ of the indices $j(1),\ldots,j(n-k)$ are contained in $\{1,\ldots,\ell\}$. Consequently, if $S_\ell(\y - \X\bb) \le \kappa_{\ell,\alpha}$, then
\begin{align*}
	\kappa_{\ell,\alpha} \
	&\ge \ G_1(\y - \X\bb) / G_\ell(\y - \X\bb) \\
	&\ge \ G_{j(1)}(\y - \X\bb) / G_{j(\ell-k)}(\y - \X\bb) \\
	&= \ G_{i(1)}(\y) / G_{i(\ell-k)}(\y) ,
\end{align*}
whence $k_{\ell,\alpha}(\y) \le \ell - 1 - (i(\ell - k) - i(1)) \le \ell - 1 - (\ell - k - 1) = k$, a contradiction to $k < k_{\ell,\alpha}(\y)$. \qed

%------------------------------------------------------------
\section{Technical details for section~\ref{subsec:Sequence}}
\label{app:Sequence}
%------------------------------------------------------------

For any $\kappa \ge 1$,
\begin{align*}
	\Pr \bigl( S_\ell(\z) \le \kappa \bigr) \
	&= \ \Pr \bigl( \tilde{\Phi}^{-1}(U_{(n)})
		\le \kappa \, \tilde{\Phi}^{-1}(U_{(n+1-\ell)}) \bigr) \\
	&= \ \Ex \Pr \bigl( U_{(n)}
		\le \tilde{\Phi}[\kappa \, \tilde{\Phi}^{-1}(U_{(n+1-\ell)})]
		\,\big|\, U_{(n+1-\ell)} \bigr) \\
	&= \ \Ex \biggl( \biggl(
		\frac{\tilde{\Phi}[\kappa \, \tilde{\Phi}^{-1}(U_{(n+1-\ell)})] - U_{(n+1-\ell)}}
		     {1 - U_{(n+1-\ell)}} \biggr)^{\ell-1} \biggr) ,
\end{align*}
where the last step follows from the fact that conditionally on $U_{(n+1-\ell)}$, the random variable $U_{(n)}$ has the same distribution as the maximum of $\ell-1$ independent random variables with uniform distirbution on $[U_{(n+1-\ell)}, 1]$. Since $U_{(n+1-\ell)}$ has distribution $\mathrm{Beta}(n+1-\ell,\ell)$, the latter expectation may be expressed as
\[
	\int_0^1
		\Bigl( \frac{\tilde{\Phi}[\kappa\,\tilde{\Phi}^{-1}(B^{-1}(u))] - B^{-1}(u)}
		            {1 - B^{-1}(u)} \Bigr)^{\ell-1} \, du ,
\]
where $B^{-1}$ is the quantile function of $\mathrm{Beta}(n+1-\ell,\ell)$. This integral can be computed numerically, and the quantile $\kappa_{\ell,\alpha}$ can be found via bisection. \qed

%------------
% References
%------------

%\bibliographystyle{ims}
%\bibliography{References}

%=============
\end{document}